\documentclass[11pt]{article}
\usepackage{anysize}\marginsize{3.5cm}{3.5cm}{1.3cm}{2cm}
\usepackage{amssymb}
\newtheorem{satz}{Satz}[section]

\newtheorem{corollary}[satz]{Corollary}
\newtheorem{definition}[satz]{Definition}
\newtheorem{example}[satz]{Example}

\newtheorem{lemma}[satz]{Lemma}
\newtheorem{proposition}[satz]{Proposition}
\newtheorem{remark}[satz]{Remark}

\newtheorem{theorem}[satz]{Theorem}
\newtheorem{conjecture}[satz]{Conjecture}

\newcommand\beginproof[1]{\trivlist\item[\hskip\labelsep{\em #1.}]}
\newcommand\proof{\beginproof{Proof}}

\def\endproof{\hspace*{\fill}\endproofsymbol\endtrivlist}

\def\endproofsymbol{\frame{\rule[0pt]{0pt}{6pt}\rule[0pt]{6pt}{0pt}}}

\renewcommand\emptyset{\varnothing}  
\renewcommand\geq{\geqslant}  
\renewcommand\leq{\leqslant}  
\renewcommand\epsilon{\varepsilon}
\renewcommand\phi{\varphi}

\renewcommand\P{\mathbb P}
\renewcommand\alpha{s}

\newcommand\be{\begin{eqnarray*}}
\newcommand\ee{\end{eqnarray*}}

\newcommand\eqnref[1]{(\ref{#1})}

\newcommand\eps{\varepsilon}

\newcommand\Z{\mathbb Z}
\newcommand\rounddown[1]{\left\lfloor#1\right\rfloor}
\newcommand\roundup[1]{\left\lceil#1\right\rceil}

\newcommand\mult{{\rm mult}}

\newcommand\newop[2]{\newcommand#1{\mathop{\rm #2}\nolimits}}
\newop\Sym{Sym}

\newcommand\call{\ell}
\newcommand\calk{k}

\begin{document}

\title{Bounds on Seshadri constants on surfaces with Picard number $1$.}
\author{Tomasz Szemberg}
\date{\today}
\maketitle
\thispagestyle{empty}

\begin{abstract}
   In this note we improve a result of Steffens \cite{Ste98} on the lower bound
   for Seshadri constants in very general points of a surface with $1$--dimensional
   N\'eron-Severi space. We also show a multi-point counterpart of such a lower bound.
\end{abstract}
\let\thefootnote\relax 

\section{Introduction}
   Seshadri constants are interesting invariants of big and nef line bundles on algebraic
   varieties. They capture the so-called local positivity of a given line bundle.
   Seshadri constants were introduced by Demailly in \cite{Dem92}. As a nice introduction
   to this circle of ideas serves \cite{PAG}, an overview of recent result is given
   in \cite{PSC}. Here we merely recall the basic definition.

\begin{definition}\label{seshcons}
   Let $X$ be a smooth projective variety, $L$ a big and nef line bundle on $X$
   and $x\in X$ a point on $X$. The number
   $$\eps(L;x):=\inf_{C\ni x}\frac{L\cdot C}{\mult_xC}$$
   is the \emph{Seshadri constant} of $L$ at $x$.
\end{definition}

   By $\eps(L;1)$ we denote the maximum
   \begin{equation}\label{eq0}
      \eps(L;1):=\max_{x\in X}\eps(L;x)
   \end{equation}
   of Seshadri constants of $L$ over all points $x\in X$. It is well known
   (see \cite[Statement 2.2.8]{PSC}) that the maximum is attained for very
   general points $x\in X$, i.e. away of a countable union of proper Zariski
   closed subsets of $X$. It is also well known (see \cite[Proposition 2.1.1]{PSC})
   that there is an upper bound
   \begin{equation}\label{upperbound}
      \eps(L; 1)\leq\sqrt[n]{L^n},
   \end{equation}
   where $n$ is the dimension of $X$.

   As for lower bounds, Steffens in \cite[Proposition 1]{Ste98} gave an interesting
   estimate on $\eps(L;1)$ in case that $X$ is a surface with Picard number $1$.
\begin{proposition}[Steffens]\label{stepro}
   Let $X$ be a smooth projective surface with Picard number $1$ and let $L$ be
   the ample generator of the N\'eron-Severi group of $X$. Then
   \begin{equation}\label{steffens}
      \eps(L;1)\geq \rounddown{\sqrt{L^2}}.
   \end{equation}
\end{proposition}

   It is clear that if $L^2$ is a square, then there is actually an equality
   $$\eps(L;1)=\sqrt{L^2}\;\;\mbox{ if }\; \sqrt{L^2}\in\Z.$$
   Our first observation is that only under these circumstances (i.e. $\sqrt{L^2}\in\Z$)
   is the bound \eqnref{steffens} sharp.

   For the rest of the paper we write $N:=L^2$ and $\alpha:=\rounddown{\sqrt{L^2}}$.

   Let as before $X$ be a smooth projective surface with Picard number $\rho(X)=1$
   and let $L$ be the ample generator of the N\'eron-Severi group.
\begin{lemma}
   If $L^2$ is not a square, then it is always
   $$\eps(L;1)>\alpha.$$
\end{lemma}
\proof
   We have by assumption that
   $\alpha < \sqrt{N}$, so that

\begin{equation}\label{eq1}
   \alpha^2+1\leq L^2.
\end{equation}
   Assume to the contrary that $\eps(L;1)=\alpha$. Then for a general point $x\in X$
   there exists a curve $C\in|pL|$, for some integer $p$, such that
\begin{equation}\label{eq2}
   pL^2=\alpha m,
\end{equation}
   where $m$ denotes as usually the multiplicity of $C$ at $x$.

   The curve $C$ cannot be smooth at $x$ because then $pL^2=\alpha$ could never be satisfied 
   (by our assumption $L^2 >1$). Hence $m\geq 2$ and we have by \cite[Theorem A]{KSS}
\begin{equation}\label{eq3}
   m(m-1)+1\leq C^2.
\end{equation}
   Combining \eqnref{eq2} and \eqnref{eq1} we get
   $$p\alpha^2<p\alpha^2+p\leq pL^2=\alpha m,$$
   which after dividing by $\alpha$ and using the fact that $\alpha$ and $p$ are integers yields
\begin{equation}\label{eq4}
   p\alpha\leq m-1.
\end{equation}
   Now, combining \eqnref{eq4} with \eqnref{eq3} we get
   $$m(m-1)+1\leq C^2=p^2L^2=p\alpha m\leq (m-1)m$$
   a contradiction.
\endproof

   With this fact established, it is natural to ask if there is a lower bound better than
   $\rounddown{\sqrt{L^2}}$ if $L^2$ is not a square. It is not obvious that such a bound exists
   because there could be a sequence of polarized surfaces $(X_n,L_n)$ with Picard number $1$, such
   that $L_n^2=N$ for all $n$ and $\lim_{n\rightarrow \infty} \eps(L_n;1)=\rounddown{\sqrt{N}}$. We show
   that this cannot happen and that there exists a lower bound on $\eps(L;1)$ improving
   that of Steffens in case $L^2$ is not a square.

\section{A new lower bound}
   We introduce some more notation. We assume that $\sqrt{N}$ is irrational and
   denote its fractional part by $\beta$,
   thus $\beta:=\sqrt{N}-\alpha>0$. We define $p_0$ as the least integer $k$ such that
   $k\cdot \beta>\frac12$, i.e.
   \begin{equation}\label{eq7}
      p_0:=\roundup{\frac{1}{2\beta}}.
   \end{equation}
   Further we set the number $m_0$ to be equal
   \begin{equation}\label{eq8}
      m_0:=\roundup{p_0\cdot\sqrt{N}}=p_0\alpha+\roundup{p_0\beta}=p_0\alpha+1.
   \end{equation}
   The following theorem is the main result of this note.
\begin{theorem}\label{main}
   Let $X$ be a smooth projective surface with Picard number $1$ and let $L$ be the
   ample generator of the N\'eron-Severi space such that $N=L^2$ is not a square. Then
   $$\eps(L;1)\geq \frac{p_0}{m_0}N.$$
\end{theorem}
\proof
   Note that $\alpha<\frac{p_0}{m_0}N<\sqrt{N}$. Indeed, as $N=(\alpha+\beta)^2$, we have
   $$\frac{p_0}{m_0}N>\frac{p_0}{p_0\alpha+1}\cdot (\alpha+2\beta)\cdot \alpha\geq \alpha.$$
   On the other hand
   $$\frac{p_0}{m_0}N=\frac{p_0\sqrt{N}}{\roundup{p_0\sqrt{N}}}\cdot\sqrt{N}<\sqrt{N}.$$
   Now, we assume to the contrary that $\eps(L;1)<\frac{p_0}{m_0} N$. Then there exists
   an integer $m$ such that for every point $x\in X$, there exists a curve $C_x$ vanishing
   at $x$ to order $\geq m$ (i.e. $\mult_xC_x\geq m$) and
   \begin{equation}\label{eq5}
      \eps(L;x)\leq\frac{L\cdot C}{m}<\frac{p_0}{m_0}N.
   \end{equation}
   Such curves $\left\{C_x\right\}$ can be chosen to form an algebraic family
   and for its arbitrary member $C$ we have
   \begin{equation}\label{eq6}
      m(m-1)+1\leq C^2
   \end{equation}
   by Theorem A in \cite{KSS}.

   On the other hand there must exist an integer $p$ such that $C\in|pL|$. The condition
   \eqnref{eq5} then translates into
   $$\frac{p}{m}<\frac{p_0}{m_0},$$
   whereas the inequality \eqnref{eq6} requires
   \begin{equation}\label{eq9}
      m(m-1)+1\leq p^2\cdot N.
   \end{equation}
   This contradicts Lemma \ref{numineq}, which we prove below.
\endproof
   We have the following numerical lemma.
\begin{lemma}\label{numineq}
   Let $N$ be a positive integer which is not a square. Let
   $$\Omega:=\left\{(p,m)\in \Z_{>0}^2:\;\; m(m-1)+1\leq Np^2\right\}$$
   and let $\eps_0:=\min_{(p,m)\in\Omega}\frac{p}{m}$. Then
   $$\eps_0=\frac{p_0}{m_0}$$
   with $p_0$ and $m_0$ defined for $N$ as in \eqnref{eq7} and \eqnref{eq8}.
\end{lemma}
\proof
   For the fixed $p$, the quotient in question is minimalized by the maximal
   integer $m$ satisfying the inequality \eqnref{eq9}. This is
   $$m_p:=\rounddown{\frac12+\sqrt{Np^2-\frac34}}.$$
   We need to show that
   $$\frac{p_0}{m_0}\leq \frac{p}{m_p}\;\; \mbox{ for all } p.$$
   We have certainly
   $$\rounddown{\frac12+p\sqrt{N}}\geq\rounddown{\frac12+\sqrt{Np^2-\frac34}}=m_p,$$
   so that it is enough to show
   \begin{equation}\label{eq10}
      \frac{p_0}{m_0}\leq \frac{p}{\rounddown{\frac12+p\sqrt{N}}}\;\; \mbox{ for all } p.
   \end{equation}
   Since
   $$p_0p\alpha+p_0\rounddown{\frac12+p\beta}=p_0\rounddown{\frac12+p\sqrt{N}}\;\;
   \mbox{ and }\; p_0m_0=p_0p\alpha+p$$
   inequality \eqnref{eq10} would follow from
   \begin{equation}\label{eq11}
      p\geq p_0\rounddown{\frac12+p\beta}.
   \end{equation}
   For $p<p_0$ the right hand side of \eqnref{eq11} is zero, since then $p\beta<\frac12$.
   For $p\geq p_0$
   we write $p=qp_0+r$ with $q\geq 1$ and $0\leq r\leq p_0-1$. In particular $r\beta<\frac12$, so that
   $$p_0\rounddown{\frac12+p\beta}=
     p_0\rounddown{\frac12+r\beta+qp_0\beta}\leq
     p_0\rounddown{qp_0\beta+1}\leq p_0q\leq p.$$
   The last but one inequality holds because $p_0\beta<1$ and $q\geq 1$.
   This verifies \eqnref{eq11} and the proof is finished.
\endproof

   The next example shows that in some situations our bound is optimal.
\begin{example}\rm
   Let $N=2d$ be an even integer such that $N+1=\ell^2$ is a square.
   A general abelian surface $X$ with polarization $L$ of type $(1,d)$
   has Picard number equal $1$, see e.g. \cite[Section 9.9]{CAV}.
   For such surfaces we know by \cite[Theorem 6.1]{Bau99} that
   $$\eps(L;1)=\frac{2d}{\ell}=\frac{1}{\sqrt{N+1}}N.$$
   On the other hand we have $p_0=1$ and $m_0=\rounddown{\sqrt{N}}+1=\ell$,
   so that
   $$\eps(L;1)=\frac{p_0}{m_0}N$$
   in that case.
\end{example}

   In general we expect however that $\eps(L;1)$ on surfaces with Picard number $1$
   is subject to a much stronger numerical restriction.
\begin{conjecture}
   Let $X$ be a smooth projective surface with Picard number $1$ and let $L$ be the
   ample generator of the N\'eron-Severi space with $N=L^2$. Then
   $$\eps(L;1)\geq \left\{
      \begin{array}{ccl}
      \sqrt{N} & \mbox{ if } & N\,\mbox{ is a square}\\
      \frac{N\calk_0}{\call_0} & \mbox{ if } & N\,\mbox{ is not a square}
      \end{array}\right.$$
   and $(\call_0,\calk_0)$ is the primitive solution of Pell's equation
   $$\call^2-N\calk^2=1.$$
\end{conjecture}
   The inequality in Theorem \ref{main} can be viewed as the next step (after Steffens) towards
   approximating $\sqrt{N}$ by continued fractions.

\section{Multi-point Seshadri constants}
   In the last paragraph we show that a lower bound of Steffens type
   can be given also for multi-point Seshadri constants. This is a variant
   of Definition \ref{seshcons} due to Xu, see \cite{Xu94}.

\begin{definition}
   Let $X$ be a smooth projective variety, $L$ a big and nef line bundle on $X$,
   $r\geq 1$ an integer and $x_1,\dots, x_r$ distinct points on $X$. The real number
   $$\eps(L; x_1,\dots, x_r):=\inf_{C\cap\left\{x_1,\dots,x_r\right\}\neq\emptyset}
     \frac{L\cdot C}{\mult_{x_1}C+\dots+\mult_{x_r}C}$$
   is the \emph{multi-point Seshadri constant} of $L$ at poins $x_1,\dots, x_r$.
\end{definition}

   The interest in these numbers comes from the fact that, at least conjecturally,
   their behavior is more predictable than that of their one-point cousins. We refere
   again to \cite[Sections 2 and 6]{PSC} for introduction to that circle of ideas.

   Similarly to \eqnref{eq0} we set
   $$\eps(L;r):=\max_{\left\{x_1,\dots,x_r\right\}\subset X} \eps(L;x_1,\dots,x_r).$$
   The following result parallels Proposition \ref{stepro}.

\begin{theorem}\label{multi}
   Let $X$ be a smooth projective surface with Picard number $1$ and let $L$ be
   the ample generator of the N\'eron-Severi group of $X$. Then
   $$\eps(L;r)\geq\rounddown{\sqrt{\frac{L^2}{r}}}.$$
\end{theorem}
\proof
   We denote $\alpha:=\rounddown{\sqrt{\frac{L^2}{r}}}$ and assume to the
   contrary that $\eps(L;r)<\alpha$.
   Then for arbitrary $x_1,\dots,x_r$ there are irreducible curves $C_{x_1,\dots,x_r}$
   such that
   \begin{equation}\label{meq1}
      \frac{L\cdot C_{x_1,\dots,x_r}}{\sum_{i=1}^r \mult_{x_i}C_{x_1,\dots,x_r}}<\alpha.
   \end{equation}
   One can choose these curves to move in an algebraic family. As the Picard number of $X$
   is $1$, there is in fact an integer $p$ such that this family is a subset of the
   linear series $pL$. If $m_1,\dots,m_r$ are positive integers such that
   $$\mult_{x_i}C_{x_1,\dots,x_r}\geq m_i\;\; \mbox{ for all }\; i=1,\dots, r,$$
   then for any member $C$ of the family we have by \cite[Lemma 1]{Xu94}
   \begin{equation}\label{meq2}
      C^2\geq m_1^2+\dots +m_{r-1}^2+m_r(m_r-1).
   \end{equation}
   We can renumber the points so that $m_r\leq m_i$ for all $i=1,\dots, r$.
   The inequality \eqnref{meq1} implies that
   $$rp\alpha^2\leq pL^2<\alpha\cdot\sum_{i=1}^r m_i.$$
   Dividing by $\alpha$ and taking into account that all involved numbers are integers
   we obtain
   \begin{equation}\label{meq3}
      rp\alpha\leq \sum_{i=1}^r m_i-1.
   \end{equation}
   On the other hand from \eqnref{meq2}, \eqnref{meq1} and \eqnref{meq3} we have
   $$\sum_{i=1}^r m_i^2-m_r\leq pL\cdot C<p\alpha\sum_{i=1}^rm_i\leq\frac1r\sum_{i=1}^rm_i\left(\sum_{i=1}^rm_i-1\right)\leq\sum_{i=1}^rm_i^2-m_r$$
   which gives the desired contradiction.
\endproof
   We have the following straightforward corollary.
\begin{corollary}
   Let $X$ and $L$ be as in Theorem \ref{multi}. Assume that the degree of $L$ is of the form
   $L^2=rd^2$ for some positive integer $d$. Then we have the equality
   $$\eps(L;r)=d.$$
\end{corollary}
\begin{remark}\rm
   The same statement was proved in \cite[Theorem 6.1.10]{PSC} for surfaces with arbitrary
   Picard number under the additional assumption that $L$ is very ample. It is expected
   that the equality
   $$\eps(L;r)=\sqrt{\frac{L^2}{r}}$$
   holds on arbitrary surfaces, provided $r$ is sufficiently large. This is a natural generalization
   of Nagata Conjecture as explained in detail in \cite{hab}.
\end{remark}

\paragraph*{Acknowledgement.}
   This work was partially supported by a MNiSW grant N N201 388834.
   I would like to thank the Max Planck Institute
   f\"ur Mathematik in Bonn, where this work has began, for warm hospitality
   and Thomas Bauer for interesting discussions. I would like also
   to thank the referee for helpful remarks.



\bigskip

   Tomasz Szemberg,
   Instytut Matematyki UP,
   Podchor\c a\.zych 2,
   PL-30-084 Krak\'ow, Poland

\nopagebreak
   \textit{E-mail address:} \texttt{szemberg@ap.krakow.pl}

\medskip
   \textit{Current address:}
   Tomasz Szemberg,
   Albert-Ludwigs-Universit\"at Freiburg,
   Mathematisches Institut,
   Eckerstra{\ss}e 1,
   D-79104 Freiburg,
   Germany.

\end{document}